\documentclass[12pt,twoside]{article}
\usepackage[english]{babel}
\usepackage[latin1]{inputenc}
\usepackage{amsmath}
\usepackage{cite}
\usepackage{amssymb,amsfonts}
\usepackage{graphicx} 
\usepackage{times,amssymb,amscd}

\newtheorem{theorem}{Theorem}[section]

\newtheorem{exmple}[theorem]{Example}
\newtheorem{defn}[theorem]{Definition}
\newtheorem{rmrk}[theorem]{Remark}

\newcommand{\bA}{\mathbf{A}}

\newcommand{\bG}{\mathbf{G}}
\newcommand{\bH}{\mathbf{H}}

\newcommand{\bL}{\mathbf{L}}
\newcommand{\bN}{\mathbf{N}}
\newcommand{\bR}{\mathbf{R}}
\newcommand{\bS}{\mathbf{S}}
\newcommand{\bV}{\mathbf{V}}
\newcommand{\bZ}{\mathbf{Z}}

\newcommand{\bY}{\mathbf{Y}}

\newcommand{\bs}{\mathbf{s}}
\newcommand{\bh}{\mathbf{h}}
\newcommand{\ba}{\mathbf{a}}
\newcommand{\bc}{\mathbf{c}}

\newcommand{\bX}{\mathbf{X}}
\newcommand{\bb}{\mathbf{b}}
\newcommand{\bu}{\mathbf{u}}

\newcommand{\BV}{\boldsymbol{V}}

\newcommand{\Bb}{\boldsymbol{b}}

\newcommand{\Bu}{\boldsymbol{u}}

\newcommand{\BF}{\boldsymbol{F}}
\newcommand{\BG}{\boldsymbol{G}}
\newcommand{\BM}{\boldsymbol{M}}

\newcommand{\cO}{\mathcal{O}}
\newcommand{\cP}{\mathcal{P}}

\newcommand{\cM}{\mathcal{M}}

\newcommand{\HYP}{\bH^3}
\newcommand{\SXR}{\bS^2\!\times\!\bR}

\newcommand{\SLR}{\widetilde{\bS\bL_2\bR}}
\newcommand{\NIL}{\mathbf{Nil}}
\newcommand{\SOL}{\mathbf{Sol}}

\begin{document}
\pagestyle{myheadings}
\markboth{\centerline{Emil Moln\'ar and Jen\H o Szirmai}}
{Infinite series of compact hyperbolic manifolds,
as possible crystal structures}
\title
{Infinite series of compact hyperbolic manifolds,
as possible crystal structures}

\author{Emil Moln\'ar and Jen\H o Szirmai \\
\normalsize Budapest University of Technology and \\
\normalsize Economics Institute of Mathematics, \\
\normalsize Department of Geometry \\
\normalsize Budapest, P. O. Box: 91, H-1521 \\
\normalsize emolnar@math.bme.hu,~szirmai@math.bme.hu
\date{\normalsize{\today}}}

\maketitle
\begin{abstract}
Previous discoveries of the first author (1984-88) on so-called hyperbolic 
football manifolds and our recent works (2016-17) on locally extremal ball 
packing and covering hyperbolic space $\HYP$ with congruent balls had led us to 
the idea that our "experience space in small size" could be of hyperbolic structure. In this paper we construct an infinite series of 
oriented hyperbolic space forms so-called cobweb (or tube) manifolds $Cw(2z, 2z, 2z)=Cw(2z)$, $3\le z$ odd, which can describe nanotubes, very probably.
\end{abstract}
%
\section{Introduction} 
\subsection{$(5,6,6)$ Archimedean solid as fundamental 
domain $\widetilde{\boldsymbol{F}}_{\BG} = \BM$ for a hyperbolic space form: Fullerene $C_{60} \sim C_{15}$}
On the base of \cite{M88} and \cite{M12} (see also \cite{CaTe09}, \cite{M84}) we shortly recall the football manifold 
$(5, 6, 6)$ in Fig.~1, as introductory example, since our infinite series 
$Cw(2z, 2z, 2z)=Cw(2z)$ $(3 \le z$ odd parameter), 
i.e. the so-called cobweb manifolds (or tube manifolds, by the later crystallographic nanotube)
will have analogous construction (Fig.~3-8). We extend also our particular result 
in \cite{MSz16}, and refer to the previous 
works cited there. This topic seems to become timely 
nowadays, also for other 3-dimensional geometries (Thurston spaces).

Our Bolyai-Lobachevsky hyperbolic space $\HYP =\cP^3\cM(\bV^4,\BV_4,\bR,\sim^, \langle~,~\rangle)$ 
will be a projective-metric space over a real vector space $\bV^4$ for points 
$X (\bX \sim c\bX,~c \in \bR \setminus\{0\})$; its dual (i.e. linear form space) $\BV_4$ 
will describe planes ($2$-planes) $u (\Bu \sim \Bu c,~ c \in \bR \setminus \{0\})$. 
The scalar product $\langle~,~\rangle$ will be specified by a 
so-called complete orthoscheme as a projective coordinate 
simplex  $b^0b^1b^2b^3 = A_0 A_1 A_2 A_3$   by $\bA_ib^j = \delta_i^j$ 
(Kronecker delta), i.e. $b^i = A_jA_kA_l, ~ \{i; j; k; l\} = \{0; 1; 2; 3\}$.

The starting Coxeter-Schl\"afli matrix will be defined by three natural parameters
$3 \le u, v, w;$ (think of $u=5, v=3, w=5$ at our football);

\begin{equation}
\begin{gathered}
(b^{ij})
=\begin{pmatrix}
1 &-\cos\frac{\pi}{u} & 0&0 \\
-\cos \frac{\pi}{u} &1& -\cos\frac{\pi}{v}&0 \\
0 & -\cos\frac{\pi}{v} &1&-\cos \frac{\pi}{w} \\
0&0&-\cos \frac{\pi}{w} &1
\end{pmatrix}=(\langle \bb^i,\bb^j \rangle),
 \end{gathered} \tag{1.1}
\end{equation}
as scalar products of basis forms $\Bb^i \in \BV_4$ $(i = 0, 1, 2, 3)$ to 
the side faces of the coordinate simplex $b^0b^1b^2b^3 = A_0 A_1 A_2 A_3$ as usual 
(see e.g. \cite{K89}, \cite{M89}, \cite{M05}, \cite{M97}, \cite{MSz}, \cite{MSz17}, \cite{V93}). 
Thus, essential face angles $(b^ib^j \angle)=\beta^{ij}$ are $\beta^{01} = \pi/u = \pi/5$, 
$\beta^{12} = \pi/v = \pi/3$, $\beta^{23} = \pi/w = \pi/5$, the others are 
$\beta^{02} = \beta^{03} = \beta^{13} = \pi/2$ (rectangle).

We assume, and this will be crucial in the following, that $u = w$, so our 
orthoscheme will be symmetric by a half-turn $\bh$. The half-turn axis $h$ joins 
the midpoints $F_{03}$ of $A_0$, $A_3$ and $F_{12}$ of $A_1$, $A_2$ (see also Fig.~2-3). 
In Fig.~1 there are $F_{12} = D_1$ and $F_{03} = D_2$ with shorter notation, the halfturn axis is denoted by $r$ there. 
At $F_{03}$ of $h$ the hexagon will be formed, orthogonal to $A_0A_3$ (other concave realization is also possible). $A_3$ (and so $A_0$ as 
well) will be the solid 
centre of $(5, 6, 6)$. $A_2$ is the face centre of a
pentagon, look also at the extended Coxeter-Schl\"afli diagrams at Fig.~1 and Fig.~2.

The last one will show our novelty at later cobweb manifolds: $\pi/u +\pi/v < \pi/2 $ and by $u = w$, 
$\pi/v +\pi/w < \pi/2$. 
Thus $A_3$ and so $A_0$ will be outer vertices, so we truncate the orthoscheme by 
their polar planes $a_3$ and $a_0$ (look also at the extended Coxeter-Schl\"afli diagram in Fig.~2), respectively, 
to obtain a compact 
domain, a so called complete orthoscheme. Algebraically, the upper minor subdeterminant sequence in (1.1) guarantees the signature 
$(+, +, +,-)$ of the scalar products (1.1), hyperbolic indeed.

Definitions of angle and distance metrics are also standard by log function of cross 
ratios, or cos and cosh functions, respectively. We mention only that the inverse 
matrix of (1.1)
\begin{equation}
(b^{ij})^{-1}=:(A_{ij})=(\langle \bA_i,\bA_j \rangle) \tag{1.2}
\end{equation}
defines the scalar product of basis vectors of $\bV^4$. These determine the distance 
metric of $\HYP$ through the scalar product $\langle \bX, \bY \rangle = X^i A_{ij} Y^j$   
of vectors $\bX = X^i \bA_i$ and $\bY = Y^j \bA_j \in \bV^4$ 
(Einstein-Schouten index conventions). Namely,
\begin{equation}
\cosh{\frac{XY}{k}}=\frac{-\langle \bX,\bY \rangle}{\sqrt{\langle \bX,\bX \rangle \langle \bY,\bY \rangle}}, ~(\langle \bX,\bX \rangle,~\langle \bY,\bY \rangle <0) \tag{1.3}
\end{equation}
is defined for the above proper points $X$, $Y$ and usually extended for other ones through complex numbers and function 
$\cos {\frac{x}{i}}= \cosh x$,  $(i=\sqrt{-1}$ is the complex imaginary unit). The natural length unit can be chosen to
$$
k=\sqrt{-\frac{1}{K}}=1, ~\text{i.e.}~K=-1
$$
is the constant negative sectional curvature.
\begin{figure}[ht]
\begin{center}
\includegraphics[width=7.5cm]{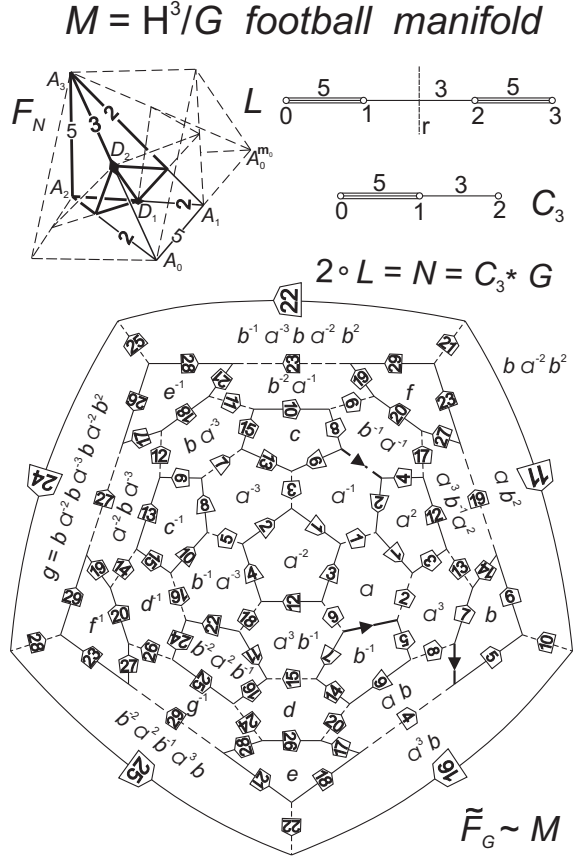}
\end{center}
\caption{The hyperbolic football manifold (fullerene), 
realized by face pairing isometries of the Archimedean solid $(5,6,6)$.}
\end{figure}
In Fig.~1 (and also in later Fig.~3-8) there are indicated the face pairing isometries
\begin{equation}
\mathbf{f}:~f^{-1}\rightarrow f, ~ ~ \mathbf{f}^{-1}:~f \rightarrow f^{-1}  \tag{1.4}
\end{equation}
of our football polyhedron $\widetilde{\BF}_{\BG}$, mapping a face $f^{-1}$ to its congruent pair $f$ by isometry $\mathbf{f},$ while 
$\widetilde{\BF}_{\BG}$ is mapped into its $f$-adjacent $\widetilde{\BF}_{\BG} ^{\mathbf{f}}$ 
(and similarly ${\mathbf{f}} ^{-1}:~  \widetilde{\BF}_{\BG} \rightarrow \widetilde{\BF}_{\BG} ^{\mathbf{f}^{-1}}$). Our main observation in \cite{M88}
was that the three starting edges, denoted by $-\!-\!\!\!\!\!\!\blacktriangleright$ in Fig.~1 had already determined 
the generator pairings
\begin{equation}
\mathbf{a}: a^{-1} \rightarrow a, ~ \mathbf{b}: b^{-1} \rightarrow b,~ \text{and their product}~ \ba\bb:~(ab)^{-1}=b^{-1}a^{-1} \rightarrow ab.    \tag{1.5}
\end{equation}
As a consequence, these induced the complete face pairing of $\widetilde{\BF}_{\BG}$, and the fundamental group ${\BG}$ of the manifold $\BM = \widetilde{\BF}_{\BG}$. We obtained $30$
edge equivalence classes (the last one is denoted by $29$), three edges in each class.

Thus, three image polyhedra join at every edge in the space tiling by $\widetilde{\BF}_{\BG}$, as the geometric presentation of the fundamental group ${\BG}$. 
Hyperbolic space $\HYP$ is just the universal cover of our football manifold $\widetilde{\BF}_{\BG} = \BM$. 
We look at Fig.~1 that the $60$ vertices of $\widetilde{\BF}_{\BG}$ fall into $15$ vertex classes under ${\BG}$, each class has $4$ vertices.
That means, $4$ footballs join at any vertex 
in the tiling. Imagine $C$ (carbon) atoms in the vertices of $\widetilde{\BF}_{\BG}$, each with $4$ valences. Therefore, $C_{15}$ would be the more logical notation instead of $C_{60}$. 
Of course, the fundamental domain $\widetilde{\BF}_{\BG}$ can contain various materials.

{\bf Summarizing}, after face pairing identification of $\widetilde{\BF}_{\BG}$ any point of it has a ball-like neighbourhood, 
or - equivalently, as well known - ${\BG}$ acts freely on $\HYP$ (without fixed point). This homogeneity seems to be advantageous for material structure.

As we discussed in \cite{M88} and \cite{M12}, ${\BG}$ is a $2$-generators group with 
presentation (defining relations to the edge classes $23$ and $26$ in Fig.~1)
\begin{equation}
\begin{gathered}
{\BG}=\{ \mathbf{a}, \bb| \mathbf{1}=(\ba^3 \bb^{-1} \ba^2 \bb^{-2}\ba^{-1})(\bb^{-2}\ba^{-1})(\bb^{-2} \ba^{2} \bb^{-1} \ba^{3}\bb)= \\ 
=(\ba^3 \bb^{-1} \ba^3 \bb)
(\bb \ba^{-2} \bb \ba^{-3}) (\bb \ba^{-2} \bb \ba^{-3} \bb \ba^{-2} \bb^{2})\}. \\
H_1(\BM)={\BG}/[{\BG},{\BG}]=\bZ_{14} \tag{1.6}
\end{gathered}
\end{equation}
is the so-called first homology group of $\BM = \widetilde{\BF}_{\BG}$ as commutator factor group of ${\BG}$ (i.e. we formally make ${\BG} \rightarrow H_1(\BM)$ commutative by 
so-called Abelianization).

The inscribed ball into $\widetilde{\BF}_{\BG}$, and so the ball packing by the tiling under group ${\BG}$, symbolizes the atomic (molecule) structure with the best known top density 
$0,77147 \dots$. Similarly, the circumscribed ball of $\widetilde{\BF}_{\BG}$ serves the best known loosest ball covering for hyperbolic space $\HYP$ \cite{MSz17} 
with density $1.36893\dots$. To this we need the generalization of volume formula of N.~I. Lobachevsky for complete orthoscheme as we cite later for information.
For other analogous ball packing and covering problems, we refer to \cite{MSzV14}, \cite{MSzV17}, \cite{Sz07-1,Sz07-2,Sz10-1,Sz11-1,Sz12,Sz13-2}, only.

We look at Fig.~1, how the half domain of orthoscheme $b^0b^1b^2b^3 = A_0 A_1 A_2 A_3$ fills out the football $(5, 6, 6) = \widetilde{\BF}_{\BG}$ by plane reflections in faces 
$b^0$, $b^1$, $b^2$  step-by-step. The well-known reflection formula for points with plane $u(\Bu)$ and its pole $U(\Bu* = \mathbf{U} = u_i b^{ij} \bA_j)$ will be
\begin{equation}
X \rightarrow Y=X^\bu,~  Y=X^\bu, ~  \bY=\bX - \frac{2(\bX\Bu)}{\langle \Bu,\Bu \rangle}\mathbf{U}. \tag{1.7}                                                              
\end{equation}
\begin{figure}[ht]
\begin{center}
\includegraphics[width=6cm]{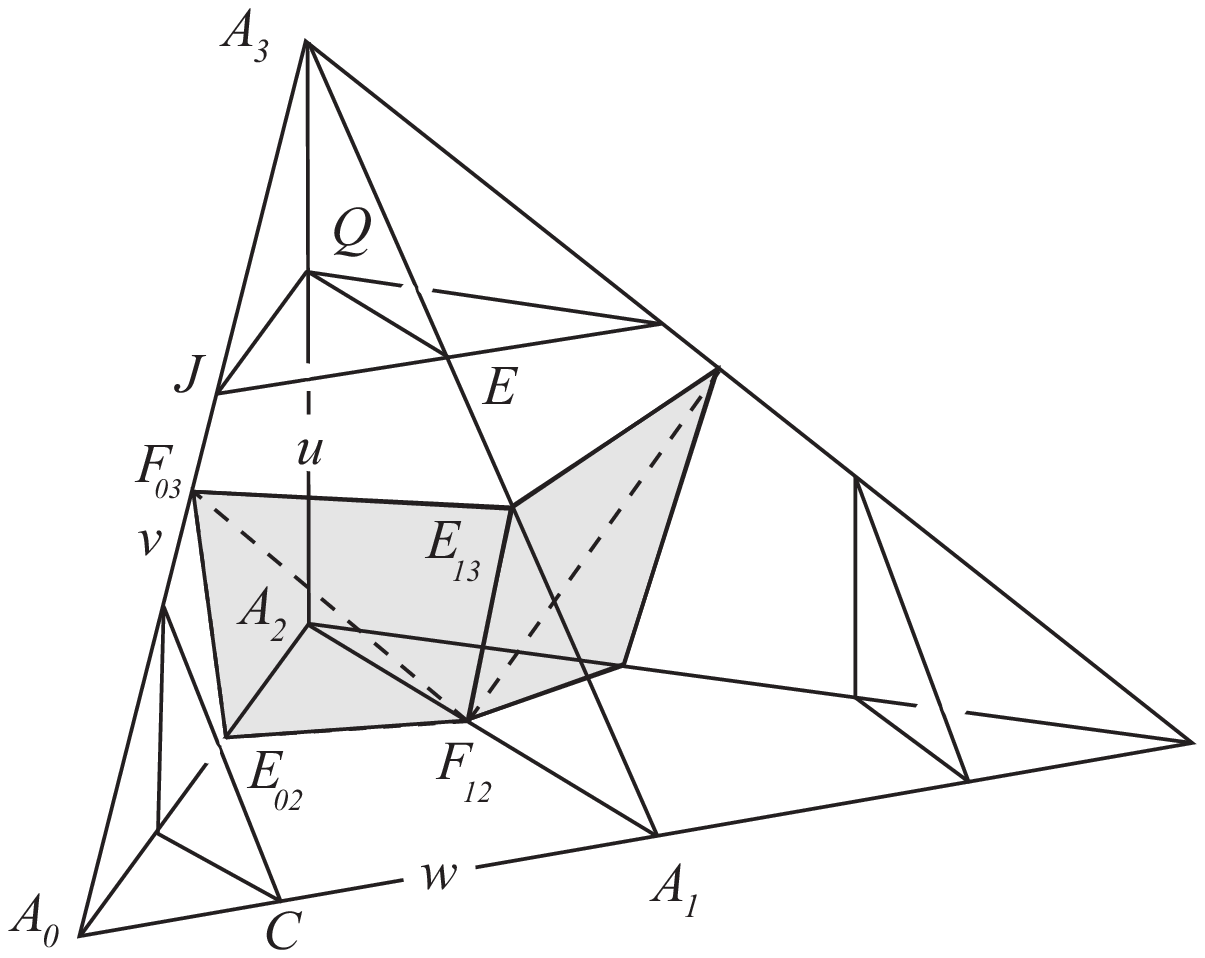} \includegraphics[width=7cm]{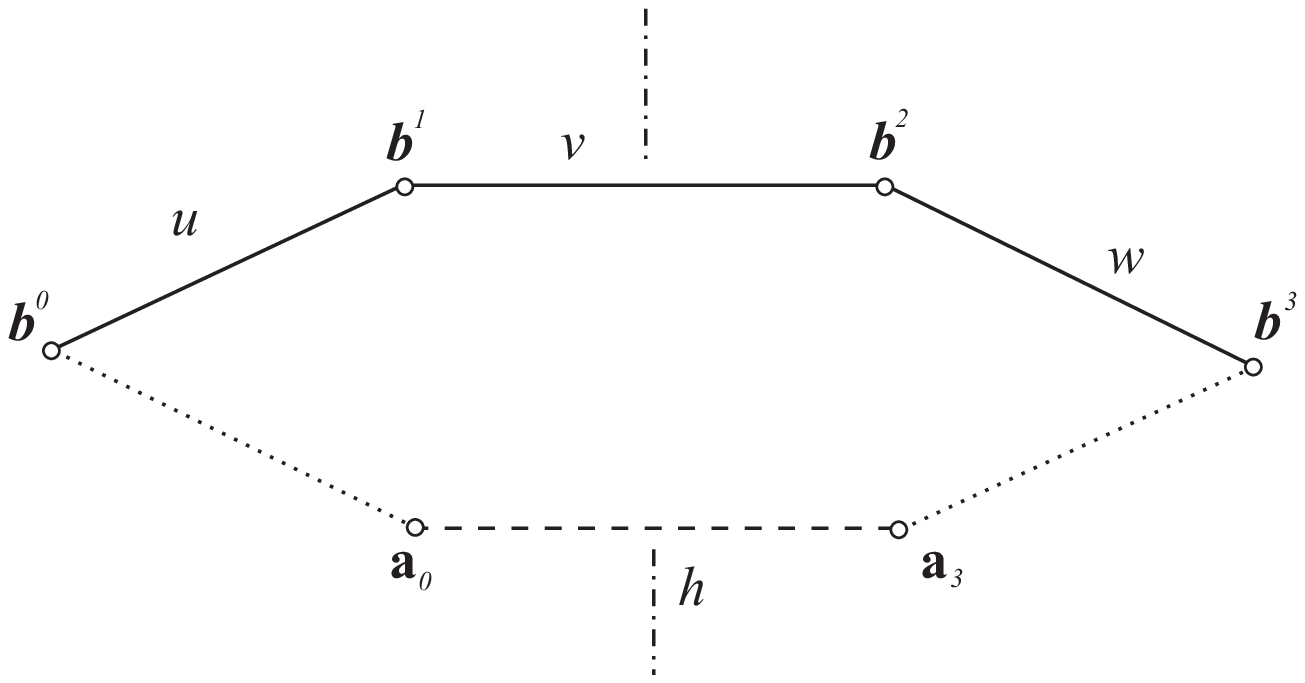}
\end{center}
\caption{The half fundamental domain of $W(u, v, w)$ for the later extended complete orthoscheme group and its extended Coxeter-Schl\"afli diagram.}
\end{figure}
Thus we get the so-called dodecahedron (or icosahedron) group with $120$ elements. 
The half-turn $\bh$, reflection $\bb^3$ and their above conjugates will map the original 
football onto the adjacent image ones. 
All these can be expressed by linear algebra as usual, and we shall use this standard 
machinery in a shorter way. 
\subsection{The cobweb (tube) manifolds $Cw(2z, 2z, 2z)$ $= Cw(2z)$, new results}
For our new cobweb manifolds (naturally named by the granddaughter Cintia of the first 
author, or we can call them {\it tube manifolds}, 
because the possible nanotube application) we shall apply 
$u = v = w = 2z$, $3 \le z$ is odd natural number, 
i.e. $b^0b^1b^2b^3 = A_0 A_1 A_2 A_3 = \mathcal{O}(2z)$ 
is an orthoscheme; complete, i.e. doubly truncated with polar planes 
$a_3$ and $a_0$ of $A_3$ and $A_0$, respectively. 

We introduce a smaller asymmetric unit $W(2z) \subset \cO(2z)$ with the half-turn axis $h$ 
and a variable halving plane through $h$. Then we choose the point 
$Q = a_3 \cap A_3A_0$ with 
its stabilizer subgroup $\bG_Q$ of order $4u = 8z$ in the extended reflection group 
$\bG$ to its fundamental domain $W(2z)$ in orthoscheme $\cO(2z)$ (Fig.~2). 

Then we reflect $W(2z)$ around $Q$ to get the cobweb polyhedron $Cw(2z)$ (Fig.~4-5) 
for a new fundamental domain of a new manifold with a new interesting fundamental 
group denoted by $\mathbf{Cw}(2z)$. 

Carbon atoms (with valence $4$, e.g. at $F_{12}$ and its $\mathbf{Cw}$-equivalent positions) can be placed very naturally in this tube-like structure.

It turns out that for $z = 4p-1$ and for $z = 4q+1$ $(1 \le p, q \in \bN )$ we get two analogous series (Fig.~4-5), each of them seems to be unique by this
manifold requirement. Our results will be formulated in Theorems 2.1, 2.2 and 3.1, 4.1 
in the next sections, respectively. But we did not succeed 
(yet?) in construction of a manifold for even $z$. These manifolds realize nanotubes in small size, very probably(!?). 
And of course, there arise new open questions.
\section{Constructions}
\subsection{Construction of cobweb (tube) manifold $Cw(6,6,6)=Cw(6)$}
By the theory, e.g in \cite{W} (cited also in works \cite{CMSpSz14}, \cite{P98}, \cite{S83}, we have to construct a fixed point free group acting in hyperbolic space $\HYP$ with the above compact fundamental domain.
In the Introduction to Fig.~1 and analogously to Fig.~2-3 we have described the extended reflection group $\bG(6,6,6)=\bG(6)$ with fundamental domain $W(6)$, as a half of 
the complete Coxeter orthoscheme $\mathcal{O}(6)$, and glued together to the cobweb polyhedron $Cw(6,6,6)$=$Cw(6)$ as Dirichlet-Voronoi (in short $D-V$) cell of the kernel point $Q$ 
by its orbit under the group $\bG(6)$. Now by Fig.~3 we shall give the face identification of $Cw(6)$, so that it will be fundamental polyhedron of the fixed-point-free group, 
denoted also by $\mathbf{Cw}(6)$, generated just by the face identifying isometries (as hyperbolic screw motions).

The complete constructioon of $Cw(6)$ has appeared in Fig.~3 with 
face pairs, signed edge triples numbered (from $1$ to $24$), signed vertex classes (by various symbols), all together $1+3 \times 3=10$ ones.

By gluing $4u=24$ domains at $Q$ (whose stabilizer subgroup $\bG_Q$ is just of order 
$|stab_Q\bG|=4u=24$) we simply "kill out" the fixed points of $\bG(6)$.
To this first, $v=u=6$ edge domains (signed by arrow $-\!\!-\!\!\!\vartriangleright$) 
is just sufficient to the former edge $F_{03}J$ of half $W(6)$ for ball-like 
neighbourhoods at points in $-\!\!-\!\!\!\vartriangleright$ edges. This can be achieved by three half-screw motions 
$\bs_1,\bs_2,\bs_3$ for the $6$ middle faces of the cobweb polyhedron $Cw(6)$, $\bs_i:s_i^{-1} \rightarrow s_i,$ $i=1,2,3$. 
The $12$ images of the former $F_{03}$ will form a vertex class $\square$, since just $24 (=4u=8z)$ domains will form the ball-like neighbourhood at these $12$ $\square$-images.

The most crucial roles are played by the former edges at the halving planes of the half orthoscheme $W(6)$ to the half-turn axis $h=F_{03}F_{12}$. 
The stabilizers of the mirror points are of order $2$ divided into two parts at $W(6)$, namely at $F_{03}E_{02}$ and at $F_{03}E_{13}$ for the odd numbered edges 
$1,3,\dots 21,23$, and of $F_{12}E_{02}$, $F_{12}E_{13}$ for the even numbered edges $2,4,\dots,22,24$, respectively. The different roles of reflection mirrors 
of $b^1$ and $b^2$, resp. $b^0$ and $b^3$, in the gluing procedure at $Q$, yield that both edge classes appear in three copies on $Cw(6)$, each class maintains 
ball-like neighborhood at every point of them.

Now comes our tricky constructions for identifying the former  half-turn faces, furthermore  the two base faces $s^{-1}$ and $s$ of $Cw(6)$ with each other (see Fig.~3). 
Two from the edge triple $1$ (to $F_{03}E_{13}$) lie on the faces $s_1^{-1}$ and $s_1$. We introduce the deciding third edge $1$ (to $F_{03}E_{12}$) 
and the orientation preserving motions $\ba_1:a_1^{-1}\rightarrow a_1$ and its inverse $\ba_1^{-1}:a_1\rightarrow a_1^{-1}$ by the mapping faces $a_1^{-1}$, $a_1$ (Fig.~3). 
This edge triple $1$ with faces $s_1^{-1},s_1,a_1^{-1},a_1$ defines a third face pairing identification $\bb_2:b_2^{-1}\rightarrow b_2$ so that
\begin{equation}
1:~\ba_1^{-1}\bs_1=\bb_2 ~ \text{holds.} \tag{2.1}
\end{equation}
Namely, three image polyhedra (between the corresponding face pairs) join each other, e.g. at the first $1$ edge in space $\HYP$ (now in combinatorial sense):
\begin{equation}
a_1^{-1}(Cw)s_1^{-1},~s_1^{\bs_1^{-1}}(Cw) ^{\bs_1^{-1}}b_2^{\bs_1^{-1}},~(b_2^{-1})^{\bb_2\bs_1^{-1}}(Cw) ^{\bb_2\bs_1^{-1}}a_1^{\bb_2\bs_1^{-1}}. \tag{2.2}
\end{equation}
Now comes again the identity polyhedron through the images
\begin{equation}
(a_1^{-1})^{\ba_1\bb_2\bs_1^{-1}}(Cw) ^{\ba_1\bb_2\bs_1^{-1}}(s_1^{-1})^{\ba_1\bb_2\bs_1^{-1}}. \tag{2.3}
\end{equation}
That means we get
\begin{equation}
\ba_1\bb_2\bs_1^{-1}=\mathbf{1} ~ \text{the identity,~equivalent~to}~ 
\Leftrightarrow \bb_2=\ba_1^{-1}\bs_1,  \notag
\end{equation}
indeed, as in (2.1).

This general combinatorial method for space filling with fundamental polyhedron, equipped by face pairing group, has been discussed algorithmically in \cite{M05}, \cite{M92}, 
\cite{MPSz06}, \cite{P92}, \cite{P98}, \cite{V93} 
in more details.
\begin{figure}[ht]
\begin{center}
\includegraphics[width=10cm]{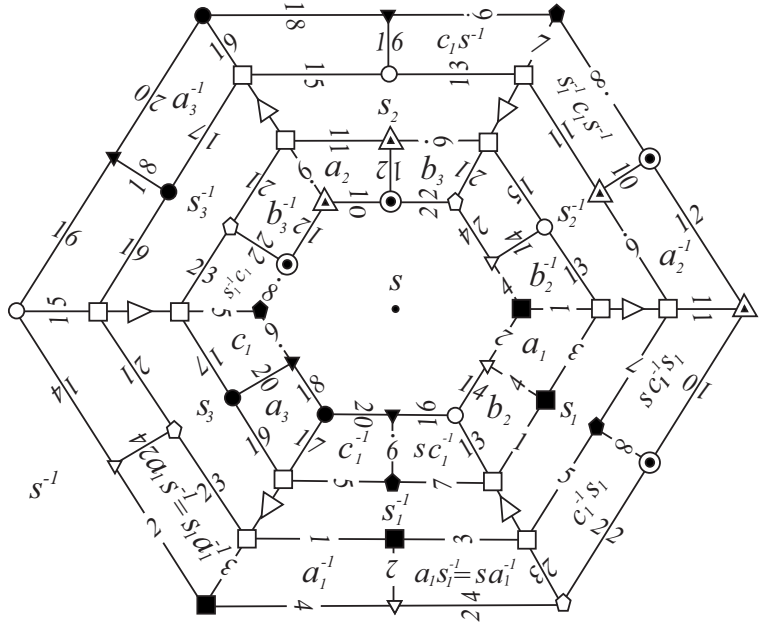}
\end{center}
\caption{The cobweb (tube) manifold $Cw(6)$ with its symbolic face pairing isometries.
Edge and vertex equivalence classes are indicated. Any point has 
a ball-like neighbourhood (nanotube).
}
\end{figure}
The first even edge class $2$ (to edges $F_{12}E_{13}$, $F_{12}E_{02}$), again in triple, just defines the most important 
identification $\bs: s^{-1} \rightarrow s$ of the base faces of
$Cw(6)$, then a new face pair:
\begin{equation}
2:~\bs\ba_1^{-1}: (sa_1^{-1})^{-1}=a_1s^{-1} \rightarrow sa_1^{-1},  \tag{2.4}
\end{equation}
with a screw motion $\bs$ through a $2\pi/3$ rotation (now). Repeating $\bs$ we get subsequent fundamental domains in $\HYP$, forming a "tube" - and finally a tiling with tubes - 
for a later material structure.
The next odd edge class $3$ defines the new face pairs with motion
\begin{equation}
3:~\ba_1 \bs_1^{-1}=\bs\ba_1^{-1} \tag{2.5}
\end{equation}
and a specific relation for this starting case $2z=6$, i.e. $u=3$, $p=1$. 
For the edge triple $4$ (to $F_{12}E_{13}$, $F_{12}E_{02}$ in Fig.2) we get a trivial relation for generators $\ba_1$ and $\bs_1$.
Our next "lucky" choice comes (from the starting triples $1$ by the cyclic $3=(4p-1)$-gonal "logical" symmetry of our cobweb polyhedron $Cw(6)$ at edge triple
$5$, where two face pairs $c_1=b_1$ and $c_1^{-1}=b_1^{-1}$ come cyclically. Then the triples $6,7,8$ follow for new pairs, transforms and a relation, respectively, as formulas 
in (2.6) show
\begin{equation}
\begin{gathered}
5:~\bs_1^{-1}\bc_1,~\text{with a formally new}~\bc_1=\bb_1,~6:~\bc_1\bs^{-1},~7:~(\bs \bc_1^{-1})\bs_1, \\ 8:~ ( \bs_1^{-1}\bc_1)\bs^{-1}(\bs\bc_1^{-1}\bs_1)=\mathbf{1}. \tag{2.6}
\end{gathered}
\end{equation}
The last relation to edge class $8$ is trivial again.
The procedure is straightforward now, and it nicely closes. To analogy of triple $1$ the edge triple $9$ defines the face pairing motion
\begin{equation}
9:~\ba_2:~ a_2^{-1} \rightarrow a_2  ~ \text{and a new motion}~ \ba_2^{-1}\bs_2=\bb_3.  \tag{2.7}
\end{equation}
The further triples $10-16$ and identifications are completely analogous. The same holds for edge triples $17-24$, starting with the face pairing motions $\ba_3$ and $\bc_1=\bb_1$
cyclically.
It turns out (see (2.5)) that our screw motions $\bs_1, \bs_2, \bs_3$  can be expressed by $\ba_1, \ba_2, \ba_3$  and $\bs$ at triples $2, 10, 18$, respectively:
\begin{equation}
\bs_1=\ba_1\bs^{-1}\ba_1,~ \bs_2=\ba_2 \bs^{-1}\ba_2,~\bs_3=\ba_3\bs^{-1}\ba_3.  \tag{2.8}
\end{equation}
\begin{figure}[ht]
\begin{center}
\includegraphics[width=10cm]{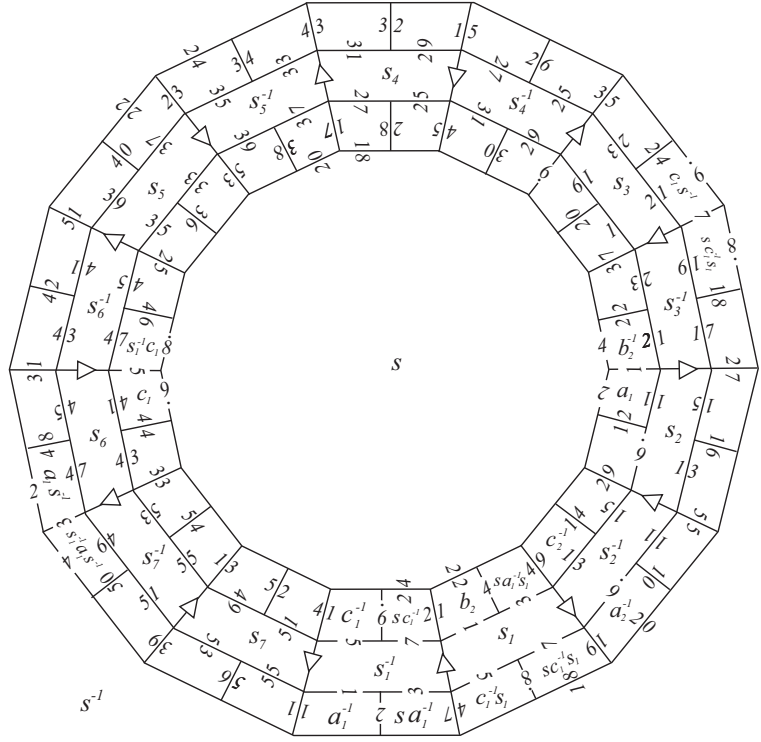}
\end{center}
\caption{The infinite series of cobweb (tube) manifolds
$Cw(2z = 8p-2)$ illustrated by $z = 7$, $p = 2$.
}
\end{figure}
The relation, to the middle edge class $-\!\!\!-\!\!\!\vartriangleright$ of six edges yields then the relation
\begin{equation}
\begin{gathered}
\mathbf{1}=(\bs_1\bs_1\bs_2\bs_2\bs_3\bs_3=(\ba_1\bs^{-1}\ba_1)^2(\ba_2\bs^{-1}\ba_2)^2(\ba_3\bs^{-1}\ba_3)^2
\end{gathered} \tag{2.9}
\end{equation}
for the fundamental group of our cobweb manifold $Cw(6)$.

But in this cyclic process, the pairing motion $\bb_1$ to edge class $5$ is not independent. Similarly 
to $1:~\bb_2 =\ba_1^{-1}\bs_1=\bs^{-1}\ba_1$, as above, we cyclically obtain $\bb_1=\ba_3^{-1}\bs_3=\bs^{-1}\ba_3$. So we get, at the edge class $7$, 
the motion $(\bs \bb_1^{-1})\bs_1=\bs(\ba_3^{-1}\bs)(\ba_1\bs^{-1}\ba_1)$ and
\begin{equation}
\begin{gathered}
\mathbf{1}=(\bs\ba_3^{-1}\bs\ba_1\bs^{-1}\ba_1)\bs\ba_2^{-1},  \tag{2.10}
\end{gathered}
\end{equation}
as well at the edge class $11$ (by $3$). At edge classes $19$ and $27$ we cyclically get
\begin{equation}
\begin{gathered}
\mathbf{1}=\bs(\ba_1^{-1}\bs)(\ba_2\bs^{-1}\ba_2)\bs\ba_3^{-1}~\text{and}~ \mathbf{1}=\bs(\ba_2^{-1}\bs)(\ba_3\bs^{-1}\ba_3)\bs\ba_1^{-1},  \tag{2.11}
\end{gathered}
\end{equation}
respectively. Eliminating $\ba_3$; we get from (2.10) and (2.11), first a $10$-letters relation, then a $18$-letters one. Keeping also the $1 \leftrightarrow 2$ 
symmetry and (2.9), we finally get three generators and three relations for the presentation of group $\mathbf{Cw}(6)$ as follow
\begin{equation}
\begin{gathered}
\mathbf{Cw}(6)=\{ \bs, \ba_1,\ba_2~| \mathbf{1}=\ba_1\ba_1\bs^{-1}\ba_1\bs\ba_2^{-1}\ba_2^{-1}\bs\ba_2^{-1}\bs^{-1}=
\\=\bs \ba_2^{-1} \bs^2 \ba_1\bs^{-1}\ba_1\bs\ba_2^{-1} \bs \ba_1^{-1} \bs \ba_2 \bs^{-1}\ba_2\bs^2\ba_1^{-1}=\\
=(\ba_1\bs^{-1}\ba_1)^2(\ba_2\bs^{-1}\ba_2)^2(\bs^{-1}\ba_2\bs^{-1}\ba_1\bs^{-1})^2\}
\end{gathered} \tag{2.12}
\end{equation}
We can introduce the first homology group 
\begin{equation}
\begin{gathered}
H_1(\mathbf{Cw}(6))=\{ \mathbf{1}=\ba_1^3 \ba_2^{-3}=\bs^6=\ba_1^6 \ba_2^{6} \bs^{2} \}
\end{gathered} \tag{2.13}
\end{equation}
of manifold $Cw(6)$, by the commutator factor group of $\mathbf{Cw}(6)$, i.e. by Abelianization.
Of course, group $\mathbf{Cw}(6)$, is a subgroup of our former $\bG(6)$ by Fig.~3. These generators $\bs, \ba_1, \ba_2$ can be expressed by the former reflections in 
$b^0$, $b^1$, $a_3$ and half-turn $h$ about $F_{03}F_{12}$ (Fig.~2).
Summarizing, we reformulate our previous partial result of \cite{MSz16} in
\begin{figure}[]
\begin{center}
\includegraphics[width=7cm]{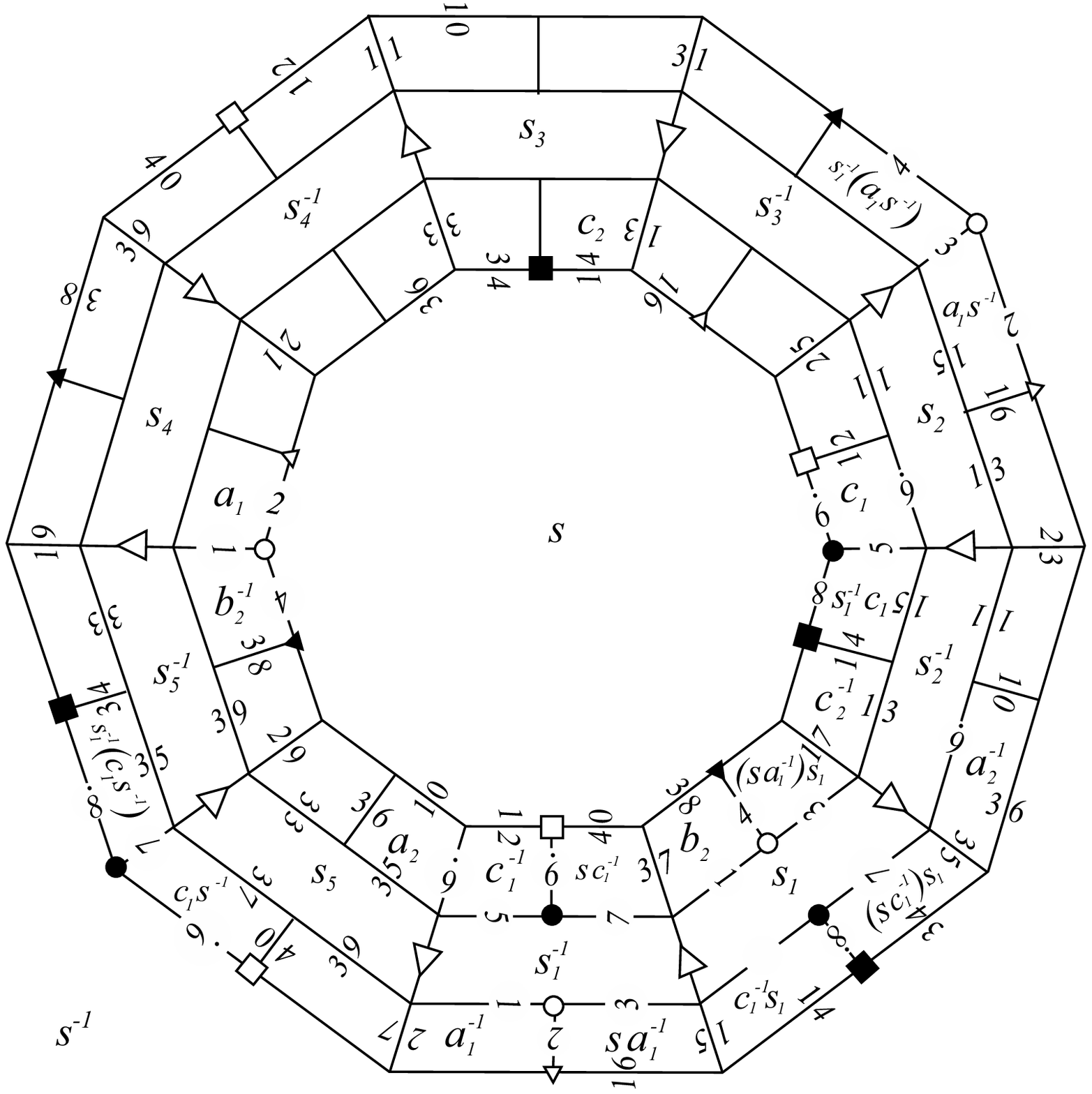} \includegraphics[width=5.5cm]{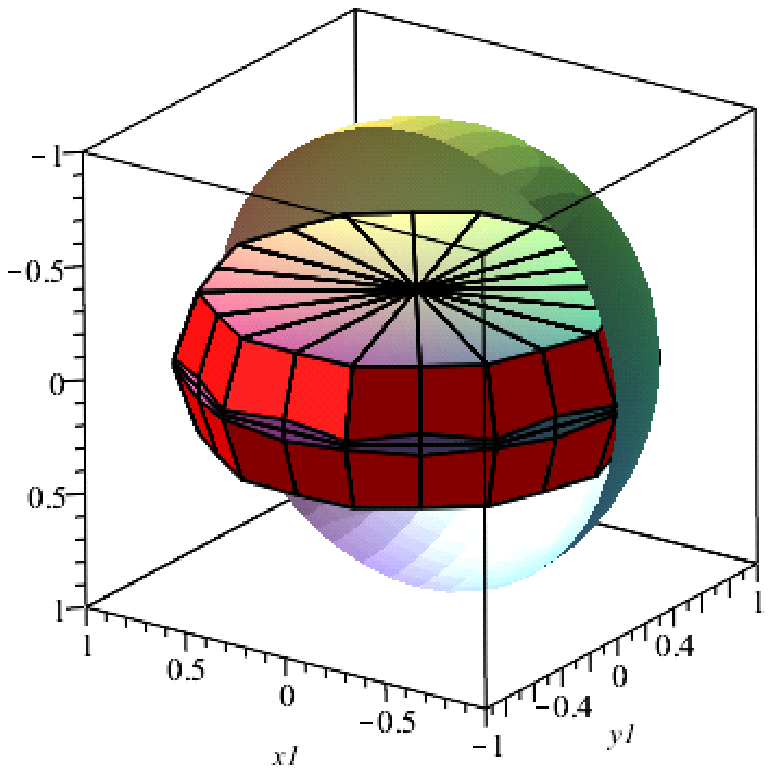}
\end{center}
\caption{Our cobweb (tube) manifolds $Cw(2z = 8q+2)$, illustrated by $z = 5$, $q = 1$.
A picture of its animation in Beltrami-Cayley-Klein model.
}
\end{figure}
\begin{theorem} 
The cobweb (tube) manifold $Cw(6)$ to cobweb polyhedron as fundamental domain has been constructed by the given face pairing identification in Fig.~3, described above. 

The fundamental group $\mathbf{Cw}(6)$ can be given by the presentation in (2.12) i.e. with three generators and three relations.  
The first homology group $H_1({Cw}(6))=H_1(\mathbf{Cw}(6))$ can be obtained by Abelianization (2.13). 

All necessary metric data of $Cw(6)$ can be computed on the base of complete orthoscheme 
$\mathcal{O}(6)$. ~ ~ $\square$
\end{theorem} 
\subsection{Construction of the cobweb (tube) manifold series $Cw(2z = 8p-2)$, $2 \le p \in \mathbf{N}$}
The above manifold $Cw(6)$ with $6 = 2z$, $z = 3 = 4p-1$, i.e. $ p = 1$ provides the analogous case $Cw(14)$ with $14 = 2z$, $z = 7 = 4p-1$, i.e. $p = 2$ (see Fig.~4). 
The screw motion $\bs$ has a rotation component $2\pi(z-1)/(2z)$, uniformly in the following. We only sketch the analogous construction, depending on parameter 
$p$, $2 \le p \in \mathbf{N}$, in general. Then we proceed obviously by induction.
The former cobweb polyhedron consists of $z=4p-1$ wedge parts. One of them is characterized by the middle face pairing (half screw motion)
\begin{equation}
\begin{gathered}
\bs_i~:~ s_i^{-1}\rightarrow s_i ~(i=1,2,\dots,z=4p-1),~\text{and so holds the relation} 
\\ \mathbf{1}=\bs_1\bs_1\bs_2\bs_2 \dots \bs_{z} \bs_{z}  \tag{2.14}
\end{gathered}
\end{equation}
to the middle edge class $-\!\!\!-\!\!\!\!\vartriangleright$. Think of $p = 2$, i.e. $7$ pairs in Fig.~4 as typical example. 
As one part contains $8$ characteristic edge classes, e.g. $1, 3, 5, 7,$ then $2, 4, 6, 8$ at the first part at 
$\bs_1: s_1^{-1} \rightarrow s_1$; we shall have all together $8\times(4p-1)=8z$ 
(now $56$) edge classes. Our sketchy analogous construction follows
\begin{equation}
\begin{gathered}
1:~\bb_2 =\ba_1^{-1}\bs_1,~2:~\bs\ba_1^{-1},~3:~\bs\ba_1^{-1}\bs_1,~5:~\bs_1^{-1}\bc_1 ~ \text{with a new}~\bc_1, \\
6:~\bc_1\bs^{-1},~ 7:~(\bs\bc_1^{-1})\bs_1,
\end{gathered} \tag{2.15}
\end{equation}
where $1,~ 3,~ 5,~ 7$ play specific roles for the later cyclic relations (indices are $\mod4p-1=7$).
\begin{equation}
\begin{gathered}
1:~\bb_{2-p} =\bc_1~ \text{one side,}~ \ba_{2-p}=\bs\ba_1^{-1}\bs_1~ \text{other side}~,~3:~\ba_{1+p}\bs^{-1}=(\bs\bc_1^{-1})\bs_1,\\
\bs\ba_{1+p}^{-1}\bs_{1+p}=\ba_2;~5:~\bc_{1+p}=\bb_2,~\bs_{1+p}^{-1}\bc_{1+p}=\bs\bc_1^{-1}; 
\\ 7:~\bc_{1-p}\bs^{-1}= \bc_1^{-1}\bs_1,~\bs_{1-p}^{-1} \bc_{1-p} \bs^{-1} 
= \bs\ba_1^{-1}.
\end{gathered} \tag{2.16}
\end{equation}
By elimination we can choose first
\begin{equation}
\begin{gathered}
\bs_{1}= \ba_1\bs^{-1} \ba_{2-p}, ~ \text{then}~ \bs_{i}= \ba_i \bs^{-1} \ba_{i+1-p},~ 
(i = 1, 2,\dots, z=4p-1),~\text{in general};\\ 
\text{then say},~ \ba_{i+1-p}= \ba_{i+3p} = \bs\ba_{i+p} \bs^{-1} \ba_{i+1} 
\bs\ba_{i+2p}^{-1} \bs, ~   (i = 1, 2,\dots, z=4p-1)
\end{gathered} \tag{2.17}
\end{equation}
\begin{equation}
\begin{gathered}
\text{besides the relation}~   \mathbf{1} = \prod_{i=1}^{z}(\ba_i \bs^{-1} \ba_{i+1-p})^2. 
\end{gathered} \tag{2.18}
\end{equation}
We briefly summarize our algorithmic result in the following
\begin{theorem}
The cobweb (tube) manifold series $Cw(2z = 8p-2)$, $2 \le p \in \bN$, for any cobweb polyhedron as $z$-cyclic fundamental domain has been algorithmically constructed above. 

The cobweb polyhedron $Cw(2z)$ is built up from the half complete orthoscheme $W(2z)$ by gluing its $8z$ copies around its vertex $Q$ as a new centre in Fig.~3-4. 

The two bases $\bs^{-1}$ and $\bs$ of this domain are paired by a screw motion $\bs:~s^{-1} \rightarrow s$ of rotational angle $2\pi(z-1)/(2z)$. 
The further face pairing is generated by $z = 4p-1$ screw motions $\ba_i:~ a_i^{-1} \rightarrow a_i$
$a_i$ given by their face pairs (indicated in Fig.~4, where $p=2$, $z=7$). 

The fundamental group is algorithmically described above 
by presentation (2.17-18). The first homology group 
$H_1({Cw}(2z))=H_1(\mathbf{Cw}(2z))$ can be obtained by Abelianization.  

All necessary metric data of $Cw(2z)$ can be computed by the complete orthoscheme $\mathcal{O}(2z)$.~ ~ $\square$
\end{theorem}
Depending on parameter $p$ we can reformulate these relations further in Section 4.
\section{Construction of cobweb (tube) manifold $Cw(10)$ and $Cw(2z = 8q + 2)$,  $1 \le q \in \bN$, in general}
In this section we construct the cobweb manifold $Cw(10; 10; 10) = Cw(10)$ in Fig.5 and analogously $Cw(8q + 2)$ for other parameters $q$ above. 
The face pairing structure of this manifold can be derived in Fig.~5. This shows the Dirichlet-Voronoi cell for group $Cw(10)$ with 
kernel point $Q$ (see Fig.~2). The above screw motion s has rotation component $2\pi(z-1)/(2z)$, throughout in the following. 
The crucial difference is that the third edge $1$ in the triple will be placed backward (opposite as above) in a $q/(4q+1)$ step on the cobweb (tube) polyhedron.
Formally, we can write just the same, but geometrically in opposite direction on the cobweb (tube):
\begin{equation}
\begin{gathered}
1:~ \bb_2 = \ba_1^{-1} \bs_1,~ 2:~ \bs\ba_1^{-1},~ 3 :~  \bs\ba_1^{-1} \bs_{1},~5:~  \bs_1^{-1} \bc_1\\
\text{with a new}~ \bc_1,~ 6:~ \bc_1 \bs^{-1},~ 7:~ (\bs\bc_1^{-1})\bs_1. 
\end{gathered} \tag{3.1}
\end{equation}
Here $1, 3, 5, 7$ play specific roles for the later cyclic relations:
\begin{equation}
\begin{gathered}
1:~ \ba_{2+q} = \bs\ba_1^{-1} \bs_{1}, ~ \bc_1= \bb_{2+q};~ 3:~ \ba_{2+q} = \bs\ba_1^{-1} \bs_{1},~ \ba_{1-q} \bs^{-1}=\bs\bc_1^{-1} s_1; \\
5:~ \bs\bc_{1+q}^{-1} = \bs_1^{-1} \bc_1, ~ \bc_1= \bb_{2+q};~ 7:~ \bs_{1+q}^{-1} \bc_{1+q} \bs^{-1}= \bs\ba_1^{-1},~ c_{1+q} \bs^{-1}= c_1^{-1} s_1.     
\end{gathered} \tag{3.2}
\end{equation}
Again we get by eliminations, first
\begin{equation}
\begin{gathered}
\bs_i = \ba_i \bs^{-1} \ba_{i+1+q},~ (i = 1, 2,\dots, z=4q + 1),~ \text{in general}; \\
\text{then say}~  \ba_{i+1+q} \bs\ba_{1-q}  \bs = \bs\ba_i^{-1} \bs^{-1} \ba_{i-2q},~ 
(i = 1, 2,\dots, z=4q + 1)~\text{and}\\
\mathbf{1} = \prod_{i=1}^{z}(\ba_i \bs^{-1} \ba_{i+1+q})^2. 
\end{gathered} \tag{3.3}
\end{equation}
provide the algorithmic presentation of the fundamental group $Cw(2z = 8q+2)$.
We summarize our algorithmic result in the following
\begin{theorem}
The cobweb (tube) manifold series $Cw(2z = 8q + 2)$, $1 \le q \in \bN$, for any cobweb polyhedron as a $z$-cyclic fundamental 
domain has been algorithmically constructed above. 
The construction is illustrated for $Cw(10)$ in Fig.~2, 5. 
The two bases $s^{-1}$ and $\bs$ of this domain are paired by a screw motion $\bs:~ s^{-1}\rightarrow \bs$ of rotational angle $2\pi(z-1)/(2z). 
The further face pairing is generated by $z = 4q+1$ screw motions \ba_i:~ a_i^{-1} \rightarrow a_i$ given by their face pairs (indicated in Fig.~5, 
where $q = 1$, $z = 5$). The fundamental group is algorithmically described above 
by presentation (3.3). The first homology group
$H_1({Cw}(2z))=H_1(\mathbf{Cw}(2z))$ can be obtained by Abelianization. 

All necessary metric data of $Cw(2z)$ can be computed by the complete orthoscheme $\mathcal{O}(2z)$.~ ~ $\square$
\end{theorem}
To this we recall only the volume formula of complete orthoscheme $\mathcal{O}(\beta^{01}, \beta^{12}, \beta^{23}$ by R.~Kellerhals to the 
Coxeter-Schl\"afli matrix (1.1) on the genial ideas of N.~I. Lobachevsky.
\begin{theorem}{\rm{(Kellerhals \cite{K89},~Lobachevsky)}} The volume of a three-dimensional hyperbolic
complete orthoscheme $\mathcal{O}(\beta^{01}, \beta^{12}, \beta^{23} \subset \mathbf{H}^3$
is expressed with the essential
angles $\beta^{01}=\frac{\pi}{u}$, $\beta^{12}=\frac{\pi}{v}$, $\beta^{23}=\frac{\pi}{w}$, $(0 \le \alpha_{ij}
\le \frac{\pi}{2})$ in the following form:

\begin{align}
&\mathrm{Vol}(\mathcal{O})=\frac{1}{4} \{ \mathcal{L}(\beta^{01}+\theta)-
\mathcal{L}(\beta^{01}-\theta)+\mathcal{L}(\frac{\pi}{2}+\beta^{12}-\theta)+ \notag \\
&+\mathcal{L}(\frac{\pi}{2}-\beta^{12}-\theta)+\mathcal{L}(\beta^{23}+\theta)-
\mathcal{L}(\beta^{23}-\theta)+2\mathcal{L}(\frac{\pi}{2}-\theta) \}, \notag
\end{align}
where $\theta \in [0,\frac{\pi}{2})$ is defined by:
$$
\tan(\theta)=\frac{\sqrt{ \cos^2{\beta^{12}}-\sin^2{\beta^{01}} \sin^2{\beta^{23}
}}} {\cos{\beta^{01}}\cos{\beta^{23}}},
$$
and where $\mathcal{L}(x):=-\int\limits_0^x \log \vert {2\sin{t}} \vert dt$ \ denotes the
Lobachevsky function (introduced by J.~Milnor in this form).
\end{theorem}
\section{New modified construction of cobweb (tube) 
manifolds $Cw(2z)$ $(3 \le z$ odd) for simplifying the fundamental 
polyhedron and the presentation of the fundamental group $\mathbf{Cw}(2z)$}
In the following we turn back to the starting complete orthoscheme $\mathcal{O}(2z)$ and its half domain $W(2z)$ and choose the halving plane 
perpendicularly to the simplex edge $A_2 A_1$. Then the point $F_{12}$ will be the common centre of $4$ previous faces of the fundamental tiling. 
Therefore, we rebuild the previous polyhedron $Cw(2z)$ so that it determines the same screw motion $\bs$ (as fixed-point-free transform). 
The faces between the previous base $\bs^{-1}$ and the stripe of the middle half-screw faces $s_i^{-1}$, $s_i$ $(i = 1, 2, \dots, z)$ will be 
glued to the previous base face $s$ just by $\bs$. Thus, the previous $s^{-1}$ and $s$ will be shifted, so that the new face $s^{-1}$ will be at the above half screw faces 
$s_i^{-1}$, $s_i$. The rotational component $2\pi(z-1)/(2z)$ makes a nice effect, it unifies the previous screw motions $\ba$, $\bb$, $\bc$ with corresponding indices.
Our Fig.~2,~6-8 show this effect as new $Cw(6)$ and $Cw(14)$, moreover $Cw(10)$ as typical examples.
\begin{figure}[ht]
\begin{center}
\includegraphics[width=14cm]{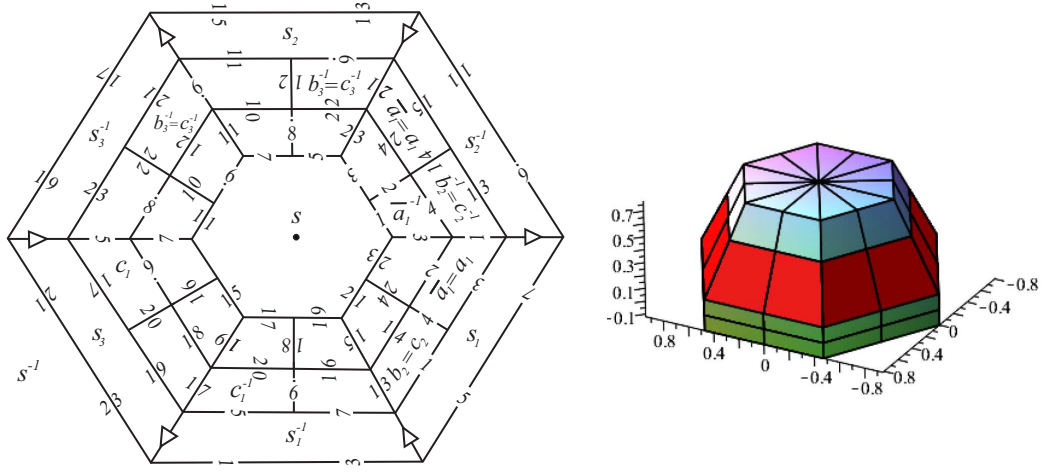}
\end{center}
\caption{$Cw(6)$ reconstructed and simplified}
\end{figure}
\begin{figure}[ht]
\begin{center}
\includegraphics[width=11cm]{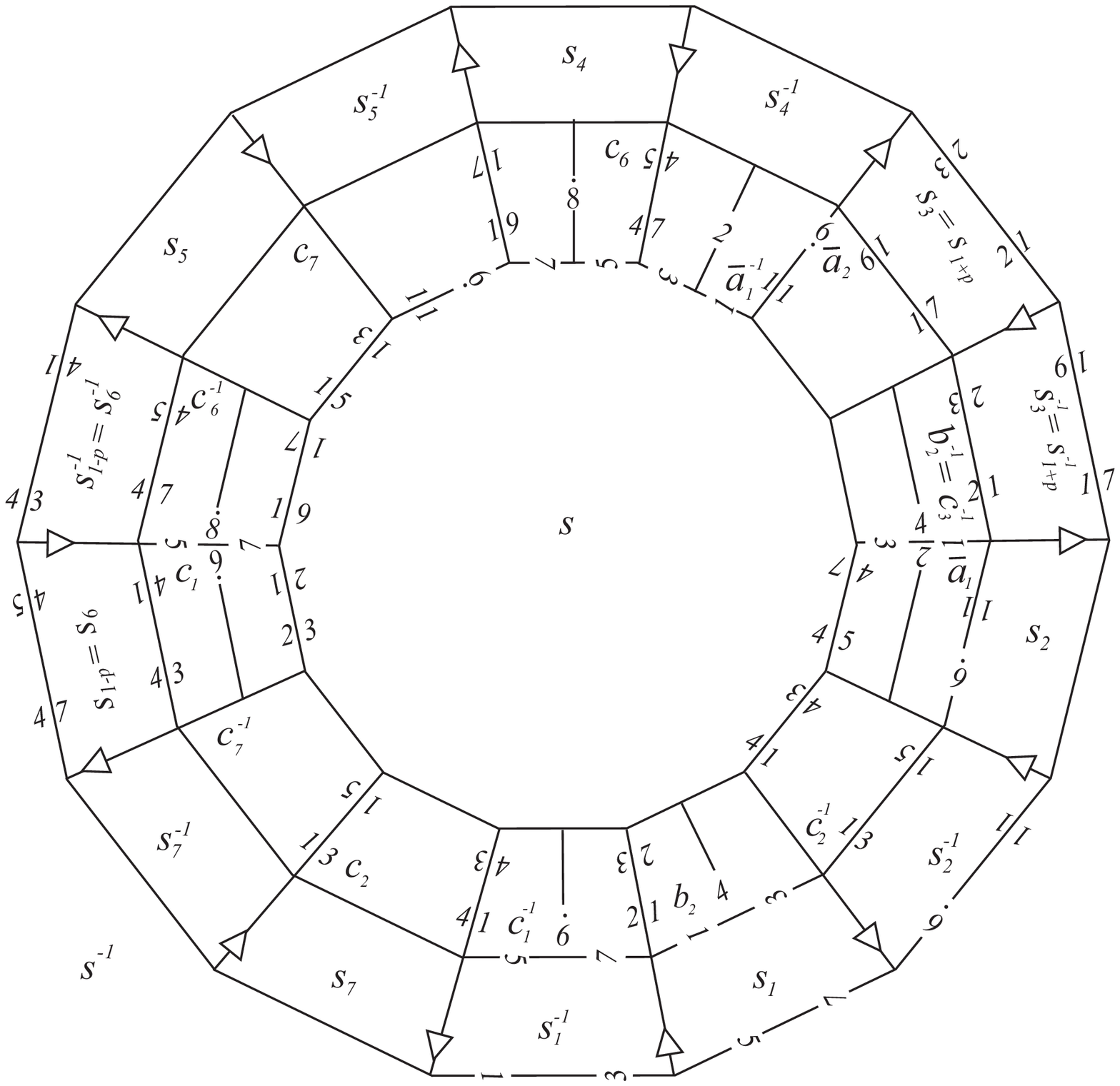}
\end{center}
\caption{$Cw(2z = 8p-2)$ illustrated by $z = 7$, $p = 2$.
Simplified reconstruction by symbolic pairing with less face pairs.}
\end{figure}
\subsection{Simplification for $Cw(2z=8p-2)$}
Let us start with the new $Cw(6)$ and $Cw(14)$ in Fig.~6-7 as $p = 1,2,$ respectively, but formulated in general, for $z = 4p-1$ (to an induction procedure). 
It turns out that we need $2z+1$ face pairs and $1 + 2z$ edge classes only. The first $-\!\!-\!\!\!\vartriangleright$ class consists of $2z$ edges to the $z$ half-turn generators
\begin{equation}
\begin{gathered}
\bs_i:~ s_i^{-1} \rightarrow s_i ~ (i=1,2, \dots, 4p-1=z)\\
\text{and so the relation}~ \mathbf{1}=\bs_1\bs_1\bs_2\bs_2\dots\bs_z\bs_z   
\end{gathered}	      \tag{4.1}
\end{equation}
in the outer strip of our cobweb, at the base $s^{-1}$. Here we find the former edge classes $1,3$ as one class, then the others, e.g. $8(i-1)+1$, $8(i-1)+3$ 
as one edge class between $s^{-1}$ and $s_i^{-1}$ $(i = 1, 2, \dots, 4p-1 = z)$. 
Please, study Fig.~7 for $Cw(14)$, i.e. $p = 2$, $z = 7$ also later on in this subsection 4.1.

The previous (Sect.~2) face pairing induced here $4$ edges in class $1,3$ and the relation now
\begin{equation}
1,3:~\bs_1b_2^{-1}\overline{\ba}_1^{-1}\bs^{-1}=\mathbf{1}~\text{equivalent to} ~ \bs_1\bc_{1+p}^{-1}\bc_{2p}^{-1}\bs^{-1}=\mathbf{1}. \tag{4.2}
\end{equation}
Similarly get to the edge classes $5,7$ the relation
\begin{equation}
\begin{gathered}
5,7:~ \bs_1^{-1} \bc_1 \bc_{1-p} \bs^{-1} = \mathbf{1},~\text{and then $\bs_1$ can be eliminated in two ways} \\
\text{to get the relation}~ \bc_{1+p}^{-1} \bc_{2p}^{-1} \bs^{-1} \bc_1 \bc_{1-p} \bs^{-1} = \mathbf{1}.    
\end{gathered}	 \tag{4.3}
\end{equation}
Or for index $i$, in general instead of index $1$, we get a simplified algorithmic presentation:
First we get for $i = 1, 2, \dots, z = 4p-1$ to the edge classes
\begin{equation}
\begin{gathered}
8(i-1)+1,~8(i-1)+3:  \bs_i \bb_{i+1}^{-1} \overline{\ba}_i^{-1} \bs^{-1} = \mathbf{1}, \\ \text{equivalent to}~ \bs_i \bc_{i+p}^{-1} \bc_{i-1+2p}^{-1} \bs^{-1}= \mathbf{1};\\
\text{then to edge classes}~8(i-1)+5,~8(i-1)+7:  \bs_i^{-1} \bc_i \bc_{i-p} \bs^{-1} =\mathbf{1}.    
\end{gathered}	 \notag
\end{equation}
Finally we obtain the desired algorithmic presentation in this $z = 4p-1~ (p = 1, 2,\dots \in \bN)$  cases:
\begin{equation}
\begin{gathered}
\bc_{i+p}^{-1} \bc_{i-1+2p}^{-1} \bs^{-1} \bc_i \bc_{i-p} \bs^{-1} = \mathbf{1}~ \text{for}~  i = 1, 2, \dots, z = 4p-1;      
\end{gathered}	 \tag{4.4}
\end{equation}
\begin{equation}
\begin{gathered}
 ~\text{and by (4.1)}~\mathbf{1} =\prod_{i=1}^{z}(\bc_i \bc_{i-p} (\bs^{-1} \bs) 
 \bc_{i-1+2p} \bc_{i+p}) \\                              
\text{so}~\mathbf{1} = \prod_{i=1}^{z}(\bc_i \bc_{i-p} \bc_{i-1+2p} \bc_{i+p}).       
\end{gathered}	 \tag{4.5}
\end{equation}
\subsection{Simplification for $Cw(2z=8q+2)$}
Let us start with the new $Cw (10)$ in Fig.~8 as $q = 1$, and formulate in general, for $z = 4q+1$ (to an inductive procedure). 
Again, we need $2z+1$ face pairs and $1 + 2z$ edge classes only. The first $-\!\!\!-\!\!\!\vartriangleright$ class consists of $2z$ edges to the $z$ half-turn generators
\begin{equation}
\begin{gathered}
\bs_i:~ s_i^{-1} \rightarrow s_i~ (i = 1, 2,\dots, 4q+1 = z)\\ \text{and to the relation}~\mathbf{1} = \bs_1 \bs_1 \bs_2 \bs_2 \dots \bs_z \bs_z       
\end{gathered}	 \tag{4.6}
\end{equation}
in the outer strip of our cobweb, at the base $s^{-1}$ as before. Here again, we find the previous edge classes $1,3$ as one class, 
then the others, e.g. $8(i-1)+1, 8(i-1)+3$ as one edge class between $s^{-1}$ and $s_i^{-1}$~ $(i = 1, 2, \dots, 4q + 1 = z)$. 
Please, study Fig.~8 for $Cw(10)$, i.e. $q = 1$, $z = 5$ also later in this subsection 4.2.

For the edge class $1,3$, then for $8(i-1)+1$, $8(i-1)+3$, in general, we get
\begin{equation}
\begin{gathered}
8(i-1)+1,~8(i-1)+3:  \bs_i \bb_{i+1}^{-1} \overline{\ba}_i^{-1} \bs^{-1} = \mathbf{1}, \\ \text{equivalent to}~ \bs_i \bc_{i-q}^{-1} \bc_{i+2q}^{-1} \bs^{-1}= \mathbf{1};\\
\text{then to edge classes}~8(i-1)+5,~8(i-1)+7:  \bs_i^{-1} \bc_i \bc_{i+q} \bs^{-1} =\mathbf{1}.    
\end{gathered}	 \tag{4.7}
\end{equation}
Finally we obtain the desired algorithmic presentation in this $z = 4q+1~ (q = 1, 2,\dots \in \bN)$  cases:
\begin{equation}
\begin{gathered}
\bc_{i-q}^{-1} \bc_{i+2q}^{-1} \bs^{-1} \bc_i \bc_{i+q} \bs^{-1} = \mathbf{1}~ 
\text{for}~  i = 1, 2, \dots, z = 4q+1;       
\end{gathered}	 \tag{4.8}
\end{equation}
\begin{equation}
\begin{gathered}
\text{and by (4.6)}~ \mathbf{1} =\prod_{i=1}^{z}(\bc_i \bc_{i+q} (\bs^{-1} \bs) \bc_{i+2q} \bc_{i-q})\\     
\text{so by simplification}~\mathbf{1} = \prod_{i=1}^{z}(\bc_i \bc_{i+q} \bc_{i+2q} \bc_{i-q}).       
\end{gathered}	 \tag{4.9}
\end{equation}
\begin{figure}[ht]
\begin{center}
\includegraphics[width=11cm]{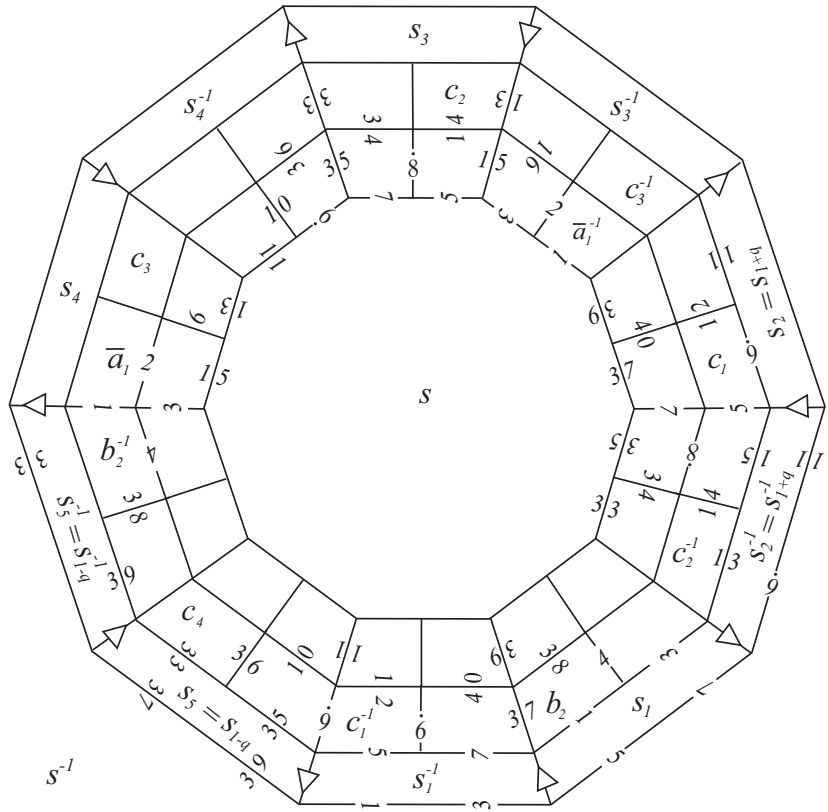}
\end{center}
\caption{Our series $Cw(2z = 8q+2)$, illustrated by $z = 5$, $q = 1$.
Symbolic reconstruction with less face pairs.}
\end{figure}
We summarize our reconstruction of cobweb (tube) manifolds in the following

\begin{theorem} The simplified cobweb (tube) manifold series $Cw(2z)$, $3 \le z$ odd 
natural number, for any reconstructed cobweb polyhedron as $z-$cyclic
fundamental domain, has been algorithmically constructed above, in Fig.~7-8, built up from the half complete orthoscheme $W(2z)$ by gluing its $8z$ copies. 

The two bases $s^{-1}$ and $s$ of this domain are paired by a screw motion $\bs:~ s^{-1} \rightarrow s$ of rotational angle $2\pi(z-1)/(2z)$. 
The further face pairing is generated by $z = 4p-1$ or $4q+1$ screw motions $\bc_i:~ c_i^{-1} \rightarrow c_{i}$ $(i = 1,2,\dots,z)$ given by their $z$ 
face pairs (indicated in Fig.~7-8, where $p = 2$, $z = 4p-1 = 7$ and $q = 1$, $z = 4q+1 = 5$; respectively).

The fundamental group is algorithmically described, above by presentation 
(4.~4-5) and (4.~8-9), respectively. 
The first homology group $H_1({Cw}(2z))=$ $H_1(\mathbf{Cw}(2z))$ can be obtained from the fundamental group by Abelianization.                                                     
All necessary metric data of $Cw(2z)$ can be computed by the complete orthoscheme $\mathcal{O}(2z)$. ~ ~ $\square$
\end{theorem}
Of course, a lot of problems arises, that we leave for other publications. For example: the above space form series seem to be 
minimal, i.e. none of them covers a smaller manifold (!?). 

For analogous problems in other Thurston geometries, we refer to 
\cite{CMSpSz14}, \cite{M05}, \cite{M97}, \cite{MSz}, \cite{MSz12}, \cite{MPSz06,MSzV14, MPSz15,  MSzV17},  \cite{S83}, \cite{Sz07-2,Sz10-1,Sz11-1,Sz12,Sz13-2}, \cite{W06}. 
We are interested in the crystallographic applications as well. Any information in this respect are kindly appreciated.

\end{document}